\documentclass[12pt,reqno,letter]{amsart}

\usepackage{amsmath,amssymb,amsthm}
\usepackage[left=2cm,top=2cm,right=2cm,bottom=2cm]{geometry}
\usepackage[normalem]{ulem}
\usepackage{mathtools}
\usepackage{enumerate}

\usepackage[a-1b]{pdfx}
\usepackage{xspace}
\usepackage{hyperref}

\newcommand{\eps}{\epsilon}   \def\t{\tau}

\newtheorem{theorem}{Theorem}[section]
\newtheorem{lemma}[theorem]{Lemma}

\theoremstyle{definition}
\newtheorem{property}[theorem]{Property}

\newcommand{\set}[1]{\left\{#1\right\}}
\def\E{\mathbb{E}}
\def\whp{{\bf whp}\xspace}

\newcommand{\beq}[2]{\begin{align}\label{#1}#2\end{align}}

\newcommand{\Pralat}{Prałat}
\newcommand{\Kaminski}{Kami\'nski}

\newcommand{\tclear}{\texttt{clear}\xspace}
\newcommand{\tstub}{\texttt{stub}\xspace}
\newcommand{\tstubn}{\texttt{stubneighbor}\xspace}
\newcommand{\tblocked}{\texttt{blocked}\xspace}

\newcommand{\stub}{{stub}\xspace} 
\newcommand{\stubs}{{stubs}\xspace} 
\newcommand{\stubv}{{stubroot}\xspace}
\newcommand{\stubvs}{{stubroots}\xspace}
\newcommand{\stubn}{{stubneighbor}\xspace}
\newcommand{\stubns}{{stubneighbors}\xspace}
\newcommand{\stubt}{{stubend}\xspace}

\newcommand{\stube}{{stubedge}\xspace}
\newcommand{\stubes}{{stubedges}\xspace}

\newcommand{\ccase}[1]{\textbf{(C\ref{#1})}}
\newcommand{\cccase}[2]{\textbf{(C\ref{#1}\textnormal{, #2})}}

\newcommand{\result}{1.85}

\definecolor{brown}{cmyk}{0, 0.72, 1, 0.45}
\definecolor{grey}{gray}{0.5}

\definecolor{lightRed}{cmyk}{0, 0.3, 0.3, 0.0}

\title{Hamilton cycles in a semi-random graph model}
\author[Alan Frieze]{Alan Frieze$^\dag$}
\address{Department of Mathematical Sciences\\Carnegie Mellon University\\Pittsburgh PA 15213}
\thanks{$\dag$ Research supported in part by NSF grant DMS1952285}
\email{alan@random.math.cmu.edu}

\author[Gregory B. Sorkin]{Gregory B. Sorkin}
\address[Gregory B. Sorkin]{Department of Mathematics,
The London School of Economics and Political Science,
Houghton Street, London WC2A 2AE, England}
\email{g.b.sorkin@lse.ac.uk}

\begin{document}

\date{5 August 2022}

\begin{abstract}
We show that \whp we can build a Hamilton cycle after at most $\result n$ rounds in a particular semi-random model. In this model, in one round,
we are given a {\em uniform random} $v\in[n]$ and then we can add an
\emph{arbitrary} edge $\set{v,w}$.
Our result improves on $2.016n$ in \cite{GMP2}.
\end{abstract}

\maketitle

\section{Introduction}
We consider the following semi-random graph model. We start with $G_0$ equal to the empty graph on vertex set $[n]$. We then obtain $G_{i+1}$ from $G_i,i\geq 0$ as follows: we are {\em presented} with a {\em uniform random} $v\in[n]$ and then we can choose to add an {\em arbitrary} edge $\set{v,w}$ to $G_i$.
This model was suggested by Peleg Michaeli
and first explored in
Ben-Eliezer, Hefetz, Kronenberg, Parczyk, Shikhelman and Stojakovi\'c \cite{BHKPSS}.
Further research on the model can be found in
Ben-Eliezer, Gishboliner, Hefetz and Krivelevich \cite{BGHK};
Gao, \Kaminski, MacRury and \Pralat \cite{GKMP};
and Gao, Macrury and \Pralat \cite{GMP1,GMP2}.
In particular \cite{GMP2} shows that \whp one can construct a Hamilton cycle in this model in at most $2.016n$ rounds.

In this short note we modify the algorithm of \cite{GMP2} and prove:
\begin{theorem}\label{th1}
In the semi-random model,
there is a strategy for constructing a Hamilton cycle in at most $\result n$ rounds.
\end{theorem}

\section{Outline and Algorithm}
Our algorithm and analysis are largely similar to those of \cite{GMP2}, so let us recapitulate the broad strokes. They maintain a large and growing path, and a set of isolated nodes. When an isolated node is presented they join it to the tail of the path.
When a path node $v$ is presented, they generate a ``\stube'' (our name, not theirs)
to a random isolated node $w$; later, if a path node $v'$ adjacent to $v$ is presented,
they generate edge $\set{v',w}$
and use it to insert
the vertex $w$ into the path between $v$ and $v'$.
These stubs are vital when the path is long and there are few isolated vertices: at that point, isolated vertices are rarely presented, while many stubs are generated. Note that by ``birthday paradox'' reasoning, only $\Theta(\sqrt n)$ \stubs are needed before there is a good chance of a neighboring vertex being presented.

As observed in \cite{GMP2},
\stubs can also be used to turn a Hamilton path into a Hamilton cycle in $o(n)$ rounds.
Assume w.l.o.g.\ that the path vertices are in sequence $1,\ldots,n$.
From each vertex $v$ presented, we generate a \stube randomly to vertex $1$ or $n$.
If $v$ had a \stube to $n$ and later $v+1$ is presented, joining $v+1$ to $1$
creates a Hamilton cycle,
using $1,\ldots,v$; $v+1,\ldots,n$;
the \stube $\set{v,n}$; and the new edge $\set{v+1,1}$.
This takes expected time $O(1/\sqrt n)$.

As in \cite{GMP2}, we maintain a large and growing path and isolated nodes,
but a key difference is that we also maintain a set of \emph{pairs}.
When an isolated node $v$ is presented, rather than joining it to the tail of the path, we join it to another isolated node $v'$ to make a {\em pair}.
When a vertex $v$ in a pair is presented, we join the pair to the tail of the path.
Stubs are used to incorporate into the path either an isolated node,
just as in \cite{GMP2}, or a pair:
If a \stube goes from $v$ to a paired vertex $w$, and a path neighbor $v'$ of $v$ is presented, joining $v'$ to the partner $w'$ of $w$ allows replacement of the edge $\set{v,v'}$ with the path $v,w,w',v'$.

The motivation for this is simple.
If an isolated vertex $v$ is presented,
the number of components decreases by 1
whether $v$ is added to the path or paired with another vertex $v'$:
in this sense, equal progress is made either way.
Ignoring the use of stubs, in the \cite{GMP2} algorithm, after $v$ is presented and joined to the path, to join $v'$ to the path we would have to wait for $v'$ to be presented.
In the paired version, after $v$ is presented and paired with $v'$,
to join the pair of them to the path we wait until \emph{either} $v$ or $v'$ is presented;
this takes half as long in expectation.

With regard to the stubs, the two versions are similar.
As just noted, the number of components (isolated vertices or pairs) needing to be incorporated into the path (including by use of stubs) is
the same either way.
The only drawback of the paired version is that
the path's growth is somewhat delayed,
so there are fewer early opportunities to create and use stubs.

\section{Algorithm and Differential Equations}
Our algorithm is largely the same as the main
``fully randomized algorithm'' of \cite{GMP2}.
We do not employ any equivalent of their initial ``degree-greedy'' phase,
although doing so would probably improve our results slightly.
Like them, we run the main algorithm until the path is
nearly but not quite complete,
so that it can be analysed by the differential equation method.
We finish up by appealing to the ``clean-up'' algorithm of \cite[Lemma 2.5]{GMP2}.

We now describe our algorithm in detail but briefly,
then present the corresponding differential equations.
As said, we maintain a \emph{path} $P$, and \emph{non-path} vertices $V$
consisting of \emph{isolated} vertices $V_1$ and \emph{paired} vertices $V_2$.

A \emph{\stube} goes from a path vertex 
we call the \emph{\stubv} or simply \emph{\stub}
to a non-path vertex 
we call a \emph{\stubt}.
We say a \stubv with $i$ \stubes has \emph{stub-degree} $i$;
we will also call it an \emph{$i$-\stubv},
and let $S_i$ be the set of such vertices.
We limit the stub-degree to at most 3, so $S=S_1\cup S_2 \cup S_3$
is the set of all \stubvs.
(There are never more than about $0.001 n$ 3-stubs,
and restricting the degree to at most 2 as in \cite{GMP2}
only increases our completion time by about $0.002 n$,
from about $1.8465 n$ to $1.8482 n$.)

Each path vertex is one of four types,
and when focussing on type we will use this font:
a \tstub;
a path-neighbor of a \stub, called a \tstubn;
a \tclear vertex, which if presented will become a \stub;
or a \tblocked vertex, which is essentially useless.

In an abuse of notation, reusing the letters for the sets to denote their cardinality,
let $P=P(t)$ denote the number of vertices on the path at time $t$, $V_1=V_1(t)$ the number of isolated vertices, $V_2=V_2(t)$ the number of vertices in pairs,
$S_i=S_i(t)$ the number of $i$-stubs (for $i \in \set{1,2,3}$),
$V(t)=V_1(t)+V_2(t)$, and $S(t)=S_1(t)+S_2(t)+S_3(t)$.
We explore the expected changes in these quantities in one round.

Our algorithm will preserve the following property.

\begin{property} \label{master} \label{pdist} \label{blocked}
Within path-distance 2 of any \tstub
there is no other \tstub nor any \tclear vertex,
and the total number of \tclear vertices is exactly $P-5S$.
%
%
\end{property}

%
%
%

We discuss this in case \ccase{C1} below.

\medskip

\subsection{Description of the Algorithm}
We list the actions taken after a (random) vertex is presented.
The description below is valid as long as there remain at least 2 isolated vertices,
and we will stop the algorithm long before that is an issue.

\begin{enumerate}[{\bf (C1)}]
\item \label{C1}
    \emph{The presented vertex $v$ is \tclear.}
    Choose a random non-path vertex $w$
    and create a \stube from $v$ to $w$,
    making $v$ a 1-\stubv.
    Change the type of $v$ from \tclear to \tstub,
    and if $v$ has path-distance 5 or more from other stubs,
    change the types of its path neighbors and second-neighbors,
    respectively, to \tstubn and \tblocked,
    i.e., \texttt{BNSNB}.
    If $v$ is at distance 4 from the next stub to the right,
    make the types from $v$ to the next stub be
    \texttt{SNBNS},
    and if distance 3, then \texttt{SNNS}.

    This changes 5 or fewer vertices from \tclear to another type.
    If fewer, then artificially change the type of additional clear
    vertices to \tblocked to make it exactly 5.
    This, and the fact that there is no \tclear nor \tstub vertex within
    path-distance 2 of $v$, preserve Property~\ref{master}.
    There is no constraint on where the artificially blocked vertices
    should be, and they need not even stay fixed from round to round.

%
%

\item \label{C2}
    \emph{The presented vertex $v$ is a \tstub with $i$ \stubns, $i=1,2$.}
    Choose a random non-path vertex $w$
    and create a \stube from $v$ to $w$,
    making $v$ an $(i+1)$-\stubv.

\item \label{C3}
    \emph{The presented vertex $v$ is a \tstubn of a stub vertex $u$.}
    By Property \ref{pdist}, $u$ is uniquely determined.
    Randomly choose one of $u$'s stubneighbors, $w$.
    Lengthen the path by removing the edge $\set{v,u}$,
    then adding the path $u,w,v$ 
    (if $w\in V_1$, making the new edge $\set{v,w}$)
    or $u,w,w',v$ 
    (if $w\in V_2$ and $w,w'$ is a pair, making the new edge $\set{v,w'}$).
    If $u$'s degree was 2 or 3, $u$ becomes an $(i-1)$-stub.

    If $u$'s degree was 1, the stub disappears:
    $u$ and its two associated \tstubn vertices become \tclear,
    as do its two associated \tblocked vertices unless they
    must remain \tblocked by proximity to some other stub.
    If necessary, change one or two other \tblocked vertices to \tclear
    to preserve Property~\ref{master}.

    The \stube $\set{v,w}$ becomes a path edge,
    and all other \stubes into $w$ (and $w'$, if relevant) are deleted.
    This results in reducing the stub-degrees of other \stubvs,
    and possibly their deletion.

\item \label{C4}
    \emph{The presented vertex $v$ is on the path,
    but \tblocked or a \stub of degree 3.}
    Do nothing.

\item \label{C5}
    \emph{The presented vertex $v$ is isolated.}
    Choose another isolated vertex $v'$ at random and make a pair $v,v'$.

\item \label{C6}
    \emph{The presented vertex $v$ is one of a pair, with some $v'$.}
    Add $v,v'$ to the tail of the path.
    As in \ccase{C3}, delete all stubs to $v$ and $v'$.
\end{enumerate}

\begin{lemma}
  \label{stubsUniform}
  In \ccase{C6}, each \stube has probability $2/V(t)$
  that its \stubt is either $v$ or $v'$,
  and these events are independent.
  In \ccase{C3} where $w$ was isolated,
  each \stube except $\set{v,w}$ has probability $1/V(t)$
  that its \stubt is $w$, and these events are independent.
  In \ccase{C3} where $w$ was paired with $w'$,
  each \stube except $\set{v,w}$ has probability $2/V(t)$ that its \stubt is either $w$ or $w'$,
  and these events are independent.
\end{lemma}

This is analogous to a claim within \cite[Lemma 2.2]{GMP2}.
For \ccase{C6} the Lemma is immediate as the stubs to $v$ and $v'$
are independent of their getting paired or joining the path.
In the \ccase{C3} cases, though,
there is a potential issue of size-biased sampling that is not explicitly
addressed in the proof in \cite{GMP2}, and so we give a proof sketch.
The issue is that $w$ being the \stubt of the chosen \stube
biases $w$ to have higher stub-degree
(e.g., $w$ could not have been selected if it had no stub edges),
suggesting that other \stubes are also more likely to have $w$ as \stubt.

\begin{proof}
Imagine that, when created, the \stubes are not revealed.
They remain, then, uniformly random between the \stubvs
(whose stub-degrees are ``known'') and non-path vertices.
Only when the \stubv $v$ is determined and one of its \stubes is chosen,
reveal (or, indeed, generate) the \stube: this determines $w$.
Only then, reveal (or generate) the other \stubes:
each is equally likely to lead to $w$ or any other non-path vertex
(including $w'$, if relevant).
So that we can apply the argument again in later rounds,
we can reveal just the \stubes incident to $w$ (and $w'$ if relevant):
after deleting them, the other \stubes remain unrevealed and uniformly random.
\end{proof}
Alternatively, one may argue from the perspective that if one sample is taken
from a population of i.i.d.\ Poisson $\lambda$ random variables,
in proportion to the variables' values,
the sampled value $X$ is not Poisson $\lambda$ (for example, it cannot be 0),
but $X-1$ is Poisson $\lambda$.

\subsection{Derivation of the Equations}
The following equations are valid as long as $P(t) \leq (1-\eps) n$
where $\eps>0$ is arbitrarily small;
anyway the differential equation method can only be applied
through such time.
The error terms below are sometimes naturally $O(1/n)$ and sometimes $O(1/V(t))$,
but with this assumption we always write them as $O(1/n)$.

\subsubsection{${\bf P(t)}$} \label{secP}
\beq{P}{
\E(P(t+1)\mid G_t)=P(t)+\frac{2V_2(t)}{n}+\frac{2S(t)}{n}\cdot \frac{V_1(t)+2V_2(t)}{V(t)} + O(1/n).
}
\begin{description}
\item[\ccase{C6}] $V_2(t)/n$ is the probability that a paired vertex is presented.
    The path length increases by 2.
\item[\ccase{C3}]
$2S(t)/n$ is the probability that a \stubn is presented.%
\footnote{Actually, a \stubn is presented with probability
$2S(t)/n+O(1/n)$, because a \stub at either end of $P$
would have only one neighbor rather than the two we are assuming.
The $O(1/n)$ correction term in \eqref{P} covers this case.
Similar correction terms apply in subsequent cases
and we will not explain the rest.}
By Property \ref{pdist} the \stubv $v$ is uniquely determined,
and one of its \stubes $\set{v,w}$ is chosen randomly.
With probability $V_1(t)/V(t)$, $w$ is isolated and the path length increases by 1;
with probability $V_2(t)/V(t)$, $w$ is paired and the path length increases by 2.
\end{description}

\subsubsection{${\bf V_1(t)}$} \label{secV1}
\beq{V1}{
\E(V_1(t+1)\mid G_t)=V_1(t)-\frac{2V_1(t)}{n}-\frac{2S(t)}{n}\cdot \frac{V_1(t)}{V(t)} + O(1/n).
}
\begin{description}
\item[\ccase{C5}] $V_1(t)/n$ is the probability that an isolated vertex is presented.
  The vertex is paired with another isolated vertex and the number of isolated vertices decreases by 2.
\item[\ccase{C3}] $2S(t)/n$ is the probability that a \stubn is presented.
  As in section \ref{secP}'s \ccase{C3}, the chosen \stubt of the \stub is isolated with probability $V_1(t)/V(t)$ and the number of isolated vertices decreases by 1.
\end{description}

\subsubsection{${\bf V_2(t)}$} \label{secV2}
\beq{V2}{
\E(V_2(t+1)\mid G_t) =V_2(t)-\frac{2V_2(t)}{n}+\frac{2V_1(t)}{n}-\frac{2S(t)}{n}\cdot \frac{2V_2(t)}{V(t)} + O(1/n).
}
\begin{description}
\item[\ccase{C6}] $V_2(t)/n$ is the probability that a paired vertex is presented.
   The path is extended using this pair, and the number of paired vertices decreases by 2.
\item[\ccase{C5}] $V_1(t)/n$ is the probability that an isolated vertex is presented.
    The vertex is paired with another isolated vertex and the number of paired vertices increases by 2.
\item[\ccase{C3}] $2S(t)/n$ is the probability that a \stubn is presented.
    The chosen \stubt of the stub is paired with probability $V_2(t)/V(t)$, and the number of paired vertices decreases by 2.
\end{description}

At this point the reader will observe that the expected change in $n=P(t)+V_1(t)+V_2(t)$ is zero, as it should be.

\subsubsection{${\bf S_1(t)}$} \label{secS1}
\beq{S1}{
\E(S_1(t+1)\mid G_t)
 &= S_1(t)+\frac{P(t)-5S(t)}{n}   - \frac{S_1(t)}{n}
 - \frac{2 S_1(t)}{n} + \frac{2 S_2(t)}{n}
 \notag \\& \hspace*{1em}
  + \frac{2S(t)}{n}
   \cdot \frac{V_1(t)+2V_2(t)}{V(t)^2}\cdot (2S_2(t) -  S_1(t))
 \notag \\& \hspace*{1em}
  + \frac{V_2(t)}{n}
   \cdot \frac{2}{V(t)}\cdot (2S_2(t) -  S_1(t))
   +O(1/n).
}
\begin{description}
\item[\ccase{C1}] $(P(t)-5S(t))/n$ is the probability that a clear vertex is presented. It becomes a 1-\stub and $S_1(t)$ increases by 1.
\item[\cccase{C2}{$i=1$}] $S_1(t)/n$ is the probability that a 1-\stub of the path is presented. It becomes a 2-\stub, and $S_1(t)$ decreases by 1.
\item[\cccase{C3}{$i=1$}] $2S_1(t)/n$ is the probability that a neighbor of a 1-\stub is presented. The \stub is used and $S_1(t)$ decreases by 1.
\item[\cccase{C3}{$i=2$}] $2S_2(t)/n$ is the probability that a neighbor of a 2-\stub is presented. The \stub is used and $S_1(t)$ increases by 1.
\item[\ccase{C3}] $2S(t)/n$ is the probability that a \stubn is presented.
    As in previous cases, the \stub $v$ is determined
    and one of its \stubes $\set{v,w}$ is chosen randomly.
    Edge $\set{v,w}$ becomes a path edge, and is no longer a \stube;
    this is captured by the previous two cases.

    All other \stubes into $w$, and its pair-partner $w'$ if any, are deleted.
    With probability $V_1(t)/V(t)$, $w$ was isolated,
    in which case by Lemma~\ref{stubsUniform}
    each \stube has probability $1/V(t)$ of having $w$ as \stubt.
    With probability $V_2(t)/V(t)$, $w$ was paired with some $w'$,
    in which case by Lemma~\ref{stubsUniform}
    each \stubt has probability $2/V(t)$ of having $w$ or $w'$ as \stubt.
    This gives the probability in the next term.

    The effect in that term is that
    each $S_2$ vertex has 2 \stubes whose potential deletion
    turns it into an $S_1$ vertex, increasing $S_1$ by 1,
    while each $S_1$ vertex has 1 \stube whose potential deletion
    turns it into a clear vertex, decreasing $S_1$ by 1.

\item[\ccase{C6}] $V_2(t)/n$ is the probability that a paired vertex is presented.
    As in the preceding case, by Lemma~\ref{stubsUniform}
    each \stube has probability $2/V(t)$ of having
    either element of the pair as \stubt and thus being deleted.
    The effect is that of the previous case.

\end{description}

\subsubsection{$\bf S_2(t)$} \label{secS2}
\beq{S2}{
\E(S_2(t+1)\mid G_t)
 &=S_2(t)+ \frac{S_1(t)}{n} - \frac{S_2(t)}{n} -\frac{2 S_2(t)}{n} + \frac{2 S_3(t)}{n}
 \notag \\& \hspace*{1em}
 +\left[
     \frac{2S(t)}{n} \cdot \frac{V_1(t)+2V_2(t)}{V(t)^2}
     + \frac{V_2(t)}{n} \cdot \frac{2}{V(t)}
  \right]
  \cdot (3S_3(t) - 2 S_2(t))
+O(1/n).
}
\begin{description}
\item[\cccase{C2}{$i=1$}] $S_1(t)/n$ is the probability that a 1-\stubv of the path is presented. It becomes a 2-\stubv, and $S_2(t)$ increases by 1.
\item[\cccase{C2}{$i=2$}] $S_2(t)/n$ is the probability that a 2-\stubv of the path is presented. It becomes a 3-\stubv, and $S_2(t)$ decreases by 1.
\item[\cccase{C3}{$i=2$}]  $2S_2(t)/n$ is the probability that a neighbor of a 2-\stubv is presented. The \stub is used and $S_2(t)$ decreases by 1.
\item[\cccase{C3}{$i=3$}] $2S_3(t)/n$ is the probability that a neighbor of a 3-\stubv is presented. The \stub is used and $S_2(t)$ increases by 1.
\item[\ccase{C3},\ccase{C6}]
    Analogous to \ccase{C3} and \ccase{C6} of section \ref{secS1},
    here combined.
\end{description}

\subsubsection{$\bf S_3(t)$} \label{secS3}
\beq{S3}{
\E(S_3(t+1)\mid G_t)
 &=S_3(t)+\frac{S_2(t)}{n}-\frac{2 S_3(t)}{n}
 \notag \\& \hspace*{1em}
 -\left[
     \frac{2S(t)}{n} \cdot \frac{V_1(t)+2V_2(t)}{V(t)^2}
     + \frac{V_2(t)}{n} \cdot \frac{2}{V(t)}
  \right]
  \cdot 3S_3(t)
     +O(1/n).
}
\begin{description}
\item[\cccase{C2}{$i=2$}] $S_2(t)/n$ is the probability that a 2-stub vertex of the path is presented. It becomes a 2-stub, and $S_3(t)$ increases by 1.
\item[\cccase{C3}{$i=3$}] $2S_3(t)/n$ is the probability that a neighbor of a 3-stub is presented. The stub is used and $S_3(t)$ decreases by 1.
\item[\ccase{C3},\ccase{C6}]
    Analogous to the corresponding case of section \ref{secS2}.
\end{description}

The equations \eqref{P} -- \eqref{S3} lead to the following differential equations in the usual way: we let $\t=t/n$ and $p(\t)=P(t)/n,v_1(\t)=V_1(t)/n$ etc. The initial conditions are $v_1(0)=1,p(0)=v_2(0)=\cdots =s_3(0)=0$.
\beq{diff}{
\begin{split}
p'&=2v_2+\frac{2s(v_1+2v_2)}{v}.\\
v_1'&=-2v_1-\frac{2sv_1}{v}.\\
v_2'&=-2v_2+2v_1-\frac{4 sv_2}{v}.\\
s_1'&=p-5s-3s_1+2s_2
 +\left[\frac{2s(v_1+2v_2)}{v^2} + \frac{2 v_2(t)}{v} \right]
      (2s_2-s_1).\\
s_2'&=s_1-3s_2+2s_3
     +\left[\frac{2s(v_1+2v_2)}{v^2} + \frac{2 v_2}{v} \right]
      (3s_3-2s_2).\\
s_3'&=s_2-2s_3
    -\left[\frac{2s(v_1+2v_2)}{v^2} + \frac{2 v_2}{v} \right]
      (3s_3).\\
\end{split}
}

\begin{figure}
  \centering
  \includegraphics[width=5.3in]{"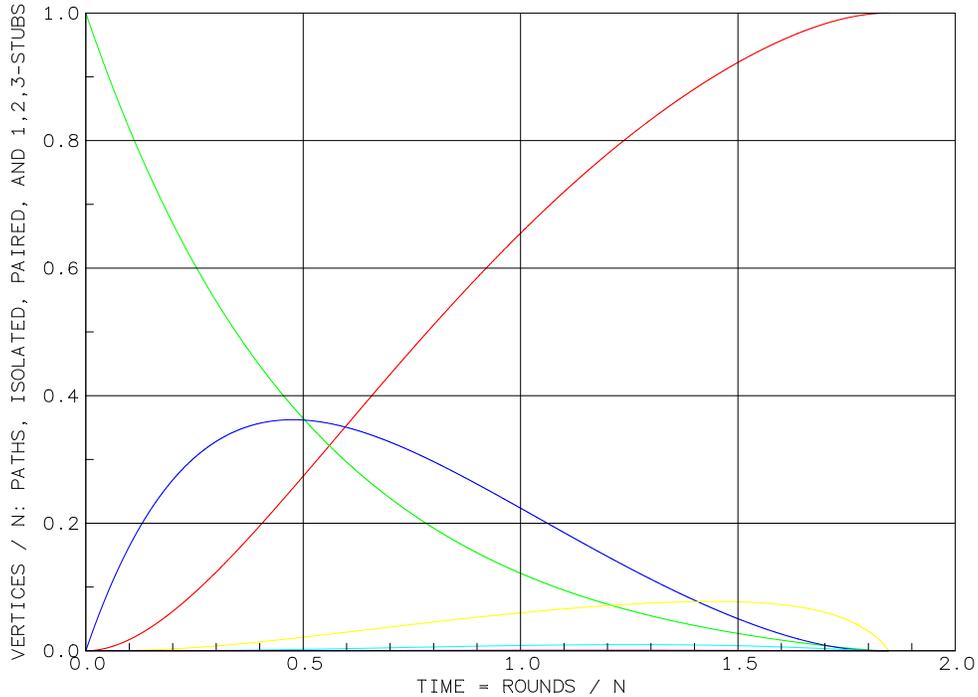"}
  \caption{The differential equations simulated,
  as a function of time (rounds divided by $n$).
  Path length $P$ (as a multiple of $n$) is shown in red,
  isolated vertices $V_1$ in green,
  paired vertices $V_2$ in blue,
  and degree 1-, 2-, and 3-stubs $S_i$ in yellow, turquoise, and magenta.
  (The figure is clearer in some PDF viewers than others.)
  }
  \label{diffEqSim}
\end{figure}

A numerical simulation of the differential equations is shown
in Figure~\ref{diffEqSim}.
It shows that $v_1(\t^*)+v_2(\t^*) \approx 0$
and $p(\t^*) \approx 1$ for $\t^*\approx \result$.
Justification of the use of the differential equation method
follows as in \cite{GMP2}.
As in \cite{GMP2}, we use the differential equation method to
analyse the algorithm until the path has length $(1-\eps)n$,
for some suitably small $\eps$.
After this we apply the \emph{clean-up} algorithm of \cite[Lemma 2.5]{GMP2}
to construct a Hamilton cycle in a further
$O(\sqrt \eps n + n^{3/4} \log^2 n)$ rounds.

\section{Concluding Remarks}

Our combining of isolated vertices into pairs leads to a substantial
speedup of the algorithm compared with \cite{GMP2},
despite our skipping their first, ``degree greedy'' phase.
We allowed for stub degrees up to 3 where \cite{GMP2} goes up only to 2,
but, observing that the number of degree-3 stubs is never more than
about $0.001 n$, this seems to have been unimportant.
Further improvements could probably be made.

First, since pairs gave a big gain,
it is natural to consider paths of 3 vertices (``triplets'')
or more.
We have not tried it, but it appears that this cannot help.
Specifically, if an isolated vertex $v$ is presented,
there would appear to be no advantage in using $v$ to extend
a ``pair'' $Q$ to a 3-vertex path,
over concatenating $v$ to the main path.
Either way, the number of components is the same.
Either way, $Q$ must eventually be brought into $P$,
either when one of its endpoints is presented
(no difference in whether $p$ is added to $Q$ or not,
as either way $Q$ has two endpoints),
or when a \stube to one of $Q$'s endpoints is used
(again, with no advantage to $v$ over $Q$'s earlier endpoint).

Other improvements, possibly challenging to analyse,
would come from choices intuitively more sensible than
the uniformly random choices made by our algorithm.

One such is to restore the ``degree greedy'' approach
from \cite{GMP2} that we discarded:
when generating \stubes, let each go to a
(random) non-path vertex \emph{of lowest \stub-degree}.

Another, when generating \stubes, is to favour paired vertices
over isolated ones.
We have some weak evidence that generating \stubes only to
non-paired vertices up to some time,
then uniformly to all non-paired vertices,
is better than generating them uniformly throughout.

Another strategy is, in the case where a 2- or 3-\stub is used,
to select a \stube
to a non-path vertex of low \stub degree,
and/or to favor an isolated vertex over a paired one
(or vice versa).

Returning to the idea of using ``triples'' as well as pairs,
potentially, small advantages could be found if, for example,
we linked $v$ with $Q$ only if $v$ were the \stubt of more \stubes
than the $Q$-endpoint it extends.

\medskip

It would be very satisfying in its own right
to better understand the natural \stub process,
where a presented vertex becomes a new \stubv unless it is already
a \stubv or \stubn, i.e., if it is at path distance at least 2 from every existing \stub.
This in contradistinction to the easier-to-analyse process
taken from \cite{GMP2} and described in \ccase{C1},
where a presented vertex becomes a stub only if it is at path distance
at least 3 from every existing \stub, and not \tblocked.
In the natural process, the number of \stubns will be between 1 and 2
times the number of \stubvs
(not 2 times, as used in \ccase{C3})
but it is not clear how to find the typical number,
nor give a good lower bound.
Presumably more \stubvs will be produced, but it is not clear how to control
the likelihood that a presented vertex will become a \stub;
indeed, nothing about the process is clear.

\section*{Acknowledgement}
We thank Zachary Hunter for a spotting two oversights in the paper's first version.

\end{document}